\documentstyle[11pt,amssymb]{article}

\begin{document}

{$\ \ \ $}
\vskip 45 pt

\centerline{\bf CLASSIFICATION OF IRREDUCIBLE REPRESENTATIONS OF}

\centerline{\bf THE q-DEFORMED ALGEBRA ${\bf U'}_{\bf q}({\bf so}_{\bf n})$}

\vskip 10 pt

\centerline{\sf A. U. KLIMYK}
\medskip

\centerline{\sc Bogolyubov Institute for Theoretical Physics}

\centerline{\sc Kiev 03143, Ukraine}

\vskip  25 pt

\begin{abstract}
A classification of finite dimensional irreducible representations of the
nonstandard $q$-deformation $U'_q({\rm so}_n)$ of the universal enveloping
algebra  $U({\rm so}(n, {\Bbb C}))$ of the Lie algebra
${\rm so}(n, {\Bbb C})$ (which does not coincides with the Drinfeld--Jimbo
quantized universal enveloping algebra $U_q({\rm so}_n)$) is given
for the case when $q$ is not a root of unity. It is shown that
such representations are exhausted by representations of the classical
and nonclassical types. Examples of the algebras $U'_q({\rm so}_3)$
and $U'_q({\rm so}_4)$ are considered in detail.
The notions of weights, highest weights, highest weight vectors are
introduced. Raising and lowering operators for irreducible finite
dimensional representations of $U'_q({\rm so}_n)$ and explicit formulas
for them are given. They depend on a weight upon which they act.
Sketch of proofs of the main assertions are given.
\end{abstract}

\section{Introduction}

Quantum orthogonal groups, quantum Lorentz groups and their
quantized universal enveloping algebras are of special interest for
modern mathematics and physics.
M. Jimbo [1] and V. Drinfeld [2] defined $q$-deformations
(quantized universal enveloping algebras)
$U_q(g)$ for all simple complex Lie algebras $g$
by means of Cartan subalgebras and root subspaces (see also [3] and [4]).
However, these approaches do not
give a satisfactory presentation of the quantized algebra
$U_q({\rm so}(n,{\Bbb C}))$ from a viewpoint of some problems in quantum
physics and mathematics.
Considering irreducible representations of the quantum
groups $SO_q(n+1)$ and $SO_q(n,1)$
we are interested in reducing them onto the quantum subgroup $SO_q(n)$.
This reduction would give an analogue of the Gel'fand--Tsetlin basis for
these representations. However, definitions of
quantized universal enveloping algebras,
mentioned above, do not allow the inclusions
$U_q({\rm so}(n+1, {\Bbb C}))\supset U_q({\rm so}(n,{\Bbb C}))$
and $U_q({\rm so}(n,1))\supset U_q({\rm so}(n))$. To be able to exploit such
reductions we have to consider $q$-deformation of the universal
enveloping algebra of the Lie algebra
${\rm so}(n+1,{\Bbb C})$ defined in terms of the generators
$I_{k,k-1}=E_{k,k-1}-E_{k-1,k}$  (where $E_{is}$ is the matrix with
elements $(E_{is})_{rt}=\delta _{ir} \delta _{st})$ rather than
by means of Cartan subalgebras and root elements.
To construct such deformations we have to deform trilinear relations
for elements $I_{k,k-1}$ instead of Serre's relations (used in the case of
quantized universal enveloping algebras of Drinfeld and
Jimbo). As a result, we obtain the associative algebra
which will be denoted as $U'_q({\rm so}_n).$

This $q$-deformation was first constructed in [5].
It permit us to construct the reductions of $U'_q({\rm so}_{n,1})$ and
$U'_q({\rm so}_{n+1})$ onto $U'_q({\rm so}_n)$.
The $q$-deformed algebra $U'_q({\rm so}_n)$ leads for $n=3$ to
the $q$-deformed algebra $U'_q({\rm so}_3)$ defined by
D. Fairlie [6]. The cyclically symmetric algebra, similar to Fairlie's
one, was also considered somewhat earlier by Odesskii [7].
The algebra
$U'_q({\rm so}_4)$ is a $q$-deformation of the
algebra $U({\rm so}(4,{\Bbb C}))$ given by means of
commutation relations between the elements $I_{ji}$, $1\le i<j\le 4$.
For the
Lie algebra ${\rm so}(4,{\Bbb C})$ we have
${\rm so}(4,{\Bbb C})={\rm so}(3,{\Bbb C})+{\rm so}(3,{\Bbb C})$,
while in the case of our $q$-deformation $U'_q({\rm so}_4)$
this is not the case (see, for example, [8]).

In the classical case, the imbedding $SO(n)\subset SU(n)$
(and its infinitesimal analogue) is of great importance for nuclear
physics and in the theory of Riemannian symmetric
spaces. It is well known
that in the framework of Drinfeld--Jimbo quantum groups and algebras
one cannot construct the corresponding embedding. The algebra
$U'_q({\rm so}_n)$ allows to define such an embedding [9],
that is, it is possible to define the embedding
$U'_q({\rm so}_n)\subset U_q({\rm sl}_n)$,
where $U_q({\rm sl}_n)$ is the Drinfeld-Jimbo quantum algebra.

As a disadvantage of the algebra $U'_q({\rm so}_n)$ we have
to mention the difficulties with Hopf algebra structure. Nevertheless,
$U'_q({\rm so}_n)$ turns out to be a coideal in
$U_q({\rm sl}_n)$ (see [9]) and this fact allows us to consider tensor
products of finite dimensional irreducible representations of
$U'_q({\rm so}_n)$ for many interesting cases (see [10]).

Finite dimensional irreducible representations of the algebra
$U'_q({\rm so}_n)$ for $q$ not a root of unity
were constructed in [5]. The formulas
of action of the generators of $U'_q({\rm so}_n)$ upon the
basis (which is a $q$-analogue of the Gel'fand--Tsetlin basis) are
given there. A proof of these formulas and some their corrections were
given in [11]. However,
finite dimensional irreducible representations described in [5] and [11]
are representations of the classical type. They are $q$-deformations of the
corresponding irreducible representations of the Lie algebra
${\rm so}_n$, that is, at $q\to 1$ they turn into representations
of ${\rm so}_n$.

If $q$ is not a root of unity,
the algebra $U'_q({\rm so}_n)$ has other classes of finite
dimensional irreducible representations which have no classical analogue.
These representations are singular at the limit $q\to 1$.
They are described in [12].
A detailed description of these
representations for the algebra $U'_q({\rm so}_3)$ is given in
[13]. A classification of irreducible $*$-representations of real forms
of the algebra $U'_q({\rm so}_3)$ is given in [14].

The aim of this paper is to give classification theorem for finite
dimensional irreducible representations of the
algebra $U'_q({\rm so}_n)$ on complex vector spaces when $q$
is not a root of unity. We show that in this case all irreducible
finite dimensional representations of $U'_q({\rm so}_n)$ are
exhausted by representations of the classical and nonclassical
types. Detailed proofs of propositions and theorems, given in this
paper, will be given separately.

Everywhere below we assume that $q$ is not a root of unity.

\section{Definition of the $q$-deformed algebra $U'_q({\rm so}_n)$}

An existence of a $q$-deformation of the
universal enveloping algebra $U({\rm so}(n,{\Bbb C}))$, different
from the Drinfeld--Jimbo
quantized universal enveloping algebra
$U_q({\rm so}_n)$, is explained by the following reason. The Lia algebra
${\rm so}(n,{\Bbb C})$ has two structures:
\bigskip

(a) The structure related to existing in ${\rm so}(n,{\Bbb C})$
a Cartan subalgebra and root elements. A quantization of this structure
leads to the Drinfeld--Jimbo
quantized universal enveloping algebra
$U_q({\rm so}_n)$.
\medskip

(b) The structure related to realization of ${\rm so}(n,{\Bbb C})$
by skew-symmetric matrices. In the Lie algebra
${\rm so}(n,{\Bbb C})$ there exists
a basis consisting of the matrices
$I_{ij}$, $i>j$, defined as $I_{ij}=E_{ij}-E_{ji}$,
where $E_{ij}$ is the matrix with entries $(E_{ij})_{rs}=
\delta _{ir}\delta _{js}$. These matrices are not root elements.
\bigskip

Using the structure (b), we may say that the universal enveloping algebra
$U({\rm so}(n,{\Bbb C}))$ is generated by the elements
$I_{ij}$, $i>j$. But in order to generate
the universal enveloping algebra
$U({\rm so}(n,{\Bbb C}))$, it is enough to take only the elements
$I_{21}$, $I_{32},\cdots ,I_{n,n-1}$. It is a minimal set of elements
necessary for generating $U({\rm so}(n,{\Bbb C}))$.
These elements satisfy the relations
$$
I^2_{i,i-1}I_{i+1,i}-2I_{i,i-1}I_{i+1,i}I_{i,i-1} +
I_{i+1,i}I^2_{i,i-1} =-I_{i+1,i},
$$   $$
I_{i,i-1}I^2_{i+1,i}-2I_{i+1,i}I_{i,i-1}I_{i+1,i} +
I^2_{i+1,i}I_{i,i-1} =-I_{i,i-1},
$$    $$
I_{i,i-1}I_{j,j-1}- I_{j,j-1}I_{i,i-1}=0\ \ \ \ {\rm for}\ \ \ \
|i-j|>1.
$$
The following theorem is true [15] for the
universal enveloping algebra $U({\rm so}(n,{\Bbb C}))$:
\bigskip

{\bf Theorem 1.}
{\it The universal enveloping algebra
$U({\rm so}(n,{\Bbb C}))$
is isomorphic to the complex associative algebra (with a unit element)
generated by the elements $I_{21}$, $I_{32},\cdots ,I_{n,n-1}$
satisfying the above relations.}
\bigskip

We make the $q$-deformation of these relations by
fulfilling the deformation of the integer 2 in these relations as
$$
2\to [2]_q:=(q^2-q^{-2})/(q-q^{-1})=q+q^{-1}.
$$
As a result, we obtain the complex unital (that is, with a unit
element) associative algebra
generated by elements $I_{21}$, $I_{32},\cdots ,I_{n,n-1}$ satisfying
the relations
$$
I^2_{i,i-1}I_{i+1,i}-(q+q^{-1})I_{i,i-1}I_{i+1,i}I_{i,i-1} +
I_{i+1,i}I^2_{i,i-1} =-I_{i+1,i}, \eqno (1)
$$    $$
I_{i,i-1}I^2_{i+1,i}-(q+q^{-1})I_{i+1,i}I_{i,i-1}I_{i+1,i} +
I^2_{i+1,i}I_{i,i-1} =-I_{i,i-1},  \eqno (2)
$$     $$
I_{i,i-1}I_{j,j-1}- I_{j,j-1}I_{i,i-1}=0\ \ \ \ {\rm for}\ \ \ \
|i-j|>1.   \eqno (3)
$$
This algebra was introduced by us in [5] and is denoted by
$U'_q({\rm so}_n)$.

The analogue of the elements $I_{ij}, i>j$, can be introduced into
$U'_q({\rm so}_n)$ (see [16]).
In order to give them we use the
notation $I_{k,k-1}\equiv I^+_{k,k-1}\equiv I^-_{k,k-1}$.
Then for $k>l+1$ we define recursively
$$
I^+_{kl}:= [I_{l+1,l},I_{k,l+1}]_{q}\equiv
q^{ 1/2}I_{l+1,l}I_{k,l+1}-
q^{- 1/2}I_{k,l+1}I_{l+1,l},   \eqno (4)
$$   $$
I^-_{kl}:= [I_{l+1,l},I_{k,l+1}]_{q^{-1}}\equiv
q^{-1/2}I_{l+1,l}I_{k,l+1}-
q^{1/2}I_{k,l+1}I_{l+1,l}.
$$
The elements $I^+_{kl}$, $k>l$, satisfy the commutation relations
$$
[I^+_{ln},I^+_{kl}]_q=I^+_{kn},\ \
[I^+_{kl},I^+_{kn}]_q=I^+_{ln},\ \
[I^+_{kn},I^+_{ln}]_q=I^+_{kl} \ \ \
{\rm for}\ \ \  k>l>n,        \eqno (5)
$$   $$
[I^+_{kl},I^+_{nr}]=0\ \ \ \ {\rm for}\ \ \
k>l>n>r\ \ {\rm and}\ \ k>n>r>l,    \eqno (6)
$$   $$
[I^+_{kl},I^+_{nr}]_q=(q-q^{-1})
(I^+_{lr}I^+_{kn}-I^+_{kr}I^+_{nl}) \ \ \ {\rm for}\ \ \
k>n>l>r.    \eqno (7)
$$
For $I^-_{kl}$, $k>l$, the commutation relations are obtained
from these relations by replacing
$I^+_{kl}$ by $I^-_{kl}$ and $q$ by $q^{-1}$.

The algebra $U'_q({\rm so}_n)$ can be defined as a unital
associative algebra generated by $I^+_{kl}$,
$1\le l<k\le n$, satisfying the relations (5)--(7). In fact, using
the relations (4) we can reduce the relations (5)--(7) to the
relations (1)--(3) for $I_{21}$, $I_{32},\cdots ,I_{n,n-1}$.

The Poincar\'e--Birkhoff--Witt theorem for the
algebra $U'_q({\rm so}_n)$ can be formulated as follows (a proof of this
theorem is given in [17]):
{\it The elements
$$
{I_{21}^+}^{m_{21}}{I_{31}^+}^{m_{31}}\cdots {I_{n1}^+}^{m_{n1}}
{I_{32}^+}^{m_{32}} {I_{42}^+}^{m_{42}} \cdots {I_{n2}^+}^{m_{n2}}
\cdots {I_{n,n-1}^+}^{m_{n,n-1}},\ \ \ \  m_{ij}=0,1,2, \cdots ,
 \eqno (8)
$$
form a basis of the algebra $U'_q({\rm so}_n)$.
This assertion is true if $I^+_{ij}$ are replaced by the
corresponding elements $I^-_{ij}$.}
\bigskip

{\bf Example 1.}
Let us consider the case of the algebra
$U'_q({\rm so}_3)$. It is generated by two elements $I_{21}$ and $I_{32}$,
satisfying the relations
$$
I_{21}^2I_{32}-(q-q^{-1})I_{21}I_{32}I_{21}+I_{32}I_{21}^2=-I_{32},  \eqno
(9)
$$    $$
I_{21}I_{32}^2-(q+q^{-1})I_{32}I_{21}I_{32}+I_{32}^2I_{21}=-I_{21}.  \eqno
(10)
$$
Introducing the element
$I^+_{31}\equiv I_{31}=q^{1/2}I_{21}I_{32}-q^{-1/2}I_{32}I_{21}$
we have for $I_{21}$, $I_{32}$, $I_{31}$ the relations
$$
[I_{21},I_{32}]_q=I_{31},\ \ \
[I_{32},I_{31}]_q=I_{21},\ \ \
[I_{31},I_{21}]_q=I_{32},   \eqno (11)
$$
where the $q$-commutator $[\cdot ,\cdot ]_q$ is defined as
$[A,B]_q=q^{1/2}AB-q^{-1/2}BA$.

Note that the algebra $U'_q({\rm so}_3)$ has a big automorphism group.
In fact, it is seen from (9) and (10) that these relations do not change
if we permute $I_{21}$ and $I_{32}$. From relations (11) we see that
the set of these relations do not change under cyclic permutation
of the elements $I_{21}$, $I_{32}$, $I_{31}$. The change of a sign
at $I_{21}$ or at $I_{32}$ also does not change the relations (9) and (10).
Generating by these automorphisms a group, we may find that they
generate the group isomorphic to the modular group $SL(2,{\Bbb Z})$.
It is why the algebra $U'_q({\rm so}_3)$ is interesting for algebraic
algebraic geometry and quantum gravity (see, for example, [18] and [19]).
\bigskip

{\bf Example 2.}
Let us consider the case of the algebra
$U'_q({\rm so}_4)$. It is generated by the elements $I_{21}$,
$I_{32}$ and $I_{43}$. We create the elements
$$
I_{31}=[I_{21},I_{32}]_q,\ \ \
I_{42}=[I_{32},I_{43}]_q,\ \ \
I_{41}=[I_{21},I_{42}]_q.    \eqno (12)
$$
Then the elements $I_{ij}$, $i>j$, satisfy the following set of relations
$$
[I_{21},I_{32}]_q=I_{31},\ \ \
[I_{32},I_{31}]_q=I_{21},\ \ \
[I_{31},I_{21}]_q=I_{32}.
$$    $$
[I_{32},I_{43}]_q=I_{42},\ \ \
[I_{43},I_{42}]_q=I_{32},\ \ \
[I_{42},I_{32}]_q=I_{43}.
$$    $$
[I_{31},I_{43}]_q=I_{41},\ \ \
[I_{43},I_{41}]_q=I_{31},\ \ \
[I_{41},I_{31}]_q=I_{43},
$$    $$
[I_{21},I_{42}]_q=I_{41},\ \ \
[I_{42},I_{41}]_q=I_{21},\ \ \
[I_{41},I_{21}]_q=I_{42},
$$   $$
[I_{21},I_{43}]=0,\ \ \ \
[I_{32},I_{41}]=0,\ \ \ \
[I_{42},I_{31}]=(q-q^{-1})(I_{21}I_{43}-I_{32}I_{41})
$$
which completely determine the algebra $U'_q({\rm so}_4)$.
At $q=1$ these relations define just the Lie algebra
${\rm so}(4,{\Bbb C})$.
Each of the sets $(I_{21},I_{32},I_{31})$,
$(I_{32},I_{43},I_{42})$, $(I_{31},I_{43},I_{41})$,
$(I_{21},I_{42},I_{41})$  determine a subalgebra isomorphic to
$U'_q({\rm so}_3)$.

The algebra $U'_q({\rm so}_4)$ is also important for quantum gravity and
algebraic geometry (see [20] and [21]). The algebra $U'_q({\rm so}_n)$
for general $n$ is also used in quantum gravity [22].
\bigskip

Let us describe the automorphism group $G$ of the algebra
$U'_q({\rm so}_n)$. It is clear from the defining relations of the
algebra $U'_q({\rm so}_n)$ that for each $i$ $(i=2,3,\cdots ,n)$
this algebra edmit an automorphism $\tau _i$ given by the formulas
\[
\tau _i: I_{j,j-1}\to I_{j,j-1},\ \ j\ne i,\ \ \ \ \ \
\tau _i: I_{i,i-1}\to -I_{i,i-1}.
\]
These automorphisms generate a group of uatomorphisms which
will be denoted by $G$. Elements of $G$ can be denoted by
$g=(\epsilon _2, \epsilon _3,\cdots ,\epsilon _n)$, where
$\epsilon _j$ runs independently the values $+1$ and $-1$.
Namely, if under action of $g$ generating elements $I_{j_1,j_1-1},
\cdots ,I_{j_s,j_s-1}$ change a sign, then in
$g=(\epsilon _2, \epsilon _3,\cdots ,\epsilon _n)$
$\epsilon _{j_1}= \cdots =\epsilon _{j_s}=-1$ and other
$\epsilon _i$ are equal to 1.
It is clear that the group $G$ has $@^{n-1}$ elements.

If $n=3$, then the group $G$ does not coincides with the group of all
automorphisms of $U'_q({\rm so}_3)$. It is not known if this assertion
is true for $n>3$.

\section{Representations of classical and nonclassical types}

The elements of the set
$$
I_{21}, I_{43},\cdots ,I_{2k,2k-1},
$$
where $n=2k$ if $n$ is even and $n=2k+1$ if $n$ is odd,
pairwise commute.
\bigskip

{\bf Proposition 1.}
(a) {\it If $T$ is a finite dimensional irreducible
representation of the algebra $U'_q({\rm so}_n)$, then the operators
$$
T(I_{21}), T(I_{43}),\cdots ,T(I_{2k,2k-1})
$$
are simulteneously diaginalizable.}

(b) {\it Possible eigenvalues of any of these operators can be as
${\rm i}[m]_q$, $m\in \frac 12 {\Bbb Z}$,
or as $[m]_+$, $m\in \frac 12 {\Bbb Z}$, $m\not\in {\Bbb Z}$,
where}
\[
[m]_q=\frac{q^m- q^{-m}}{q-q^{-1}},\ \ \ \
[m]_+=\frac{q^m+ q^{-m}}{q-q^{-1}}.
\]
\medskip

Eigenvalues of the form ${\rm i}[m]_q$ are called {\it eigenvalues of the
classical type}.
Eigenvalues of the form $[m]_+$ are called {\it eigenvalues of the
nonclassical type}.

The following proposition is important for construction of weight
theory for finite dimensional representations of $U'_q({\rm so}_n)$.
\bigskip

{\bf Proposition 2.}
{\it Let $T$ be a finite dimensional irreducible
representation of $U'_q({\rm so}_n)$. Then}

(a) {\it Eigenvalues of any operator $T(I_{2i,2i-1})$ are all of the
classical type or all of the nonclassical type.}

(b) {\it Moreover, all operators $T(I_{2i,2i-1})$, $i=1,2,\cdots ,k$, have
eigenvalues of the same type.}
\bigskip

This proposition is proved by restricting the representation $T$ to
the subalgebras $U'_q({\rm so}_4)$ generated by the elements
$I_{j,j-1}, I_{j+1,j}, I_{j+2,j+1}$, $j=2,3,\cdots ,n-2$ and using the
results of the paper [8].
\bigskip

{\bf Definition 1.}
A finite dimensional irreducible representation $T$
of the algebra $U'_q({\rm so}_n)$ is called of {\it classical}
({\it nonclassical}) {\it type} if the operators
$T(I_{2i,2i-1})$, $i=1,2,\cdots ,k$ have eigenvalues of the
classical (of the nonclassical) type.
\bigskip

{\bf Proposition 3.}
{\it Let $T$ be a finite dimensional irreducible
representation of $U'_q({\rm so}_n)$ of the classical (nonclassical) type.
Then a restriction of $T$ to the subalgebra $U'_q({\rm so}_{n-1})$
decomposes into a direct sum of irreducible representations of
this subalgebra belonging to the same type.}

\section {Weights of representations}

In this section we construct a $q$-analogue of weights for
finite dimensional irreducible representations of the algebra
$U'_q({\rm so}_n)$. Note that this algebra has no elements which can be
treated as root elements (similar to root elements of semisimple Lie
algebras or quantized universal enveloping algebras of Drinfeld and Jimbo).
For this reason, we have not a weight theory for finite dimensional
representations of $U'_q({\rm so}_n)$ similar to that for semisimple Lie
algebras. However, we can construct the theory which can
replace the weight theory of representations of semisimple Lie algebras.
\bigskip

{\bf Definition 2.}
Let $T$ be a finite dimensional representation of the
algebra $U'_q({\rm so}_n)$. Eigenvectors ${\bf v}$ of operators
$T(I_{2j,2j-1})$,
$j=1,2,\cdots ,k$, are called {\it weight vectors} of the representation
$T$. If $T(I_{2j,2j-1}){\bf v}=m_j {\bf v}$, then the
the set of numbers ${\bf m}=(m_1,m_2,\cdots ,m_{k})$, where $n=2k+1$
or $n=2k$, is called a weight of the vector ${\bf v}$.
\bigskip

The set of all weights of an irreducible representation $T$ of
$U'_q({\rm so}_n)$ is called a {\it weight diagram} of the representation
$T$.
\bigskip

{\bf Proposition 4.}
{\it A weight diagram of a finite dimensional irreducible representation $T$
of the classical type is invariant with respect to the Weyl group $W$ of the
Lie algebra ${\rm so}(n,{\Bbb C})$.}
\bigskip

This proposition is proved by restriction of the representation $T$
to the subalgebras $U'_q({\rm so}_3)$ generated by the pairs of
generators $I_{2j,2j-1}, I_{2j+1,2j}$, $j=1,2,\cdots$, and using the
results of the paper [23].

Note that
a weight diagram of a finite dimensional irreducible representation
of the nonclassical type is not invariant with respect to the Weyl
group $W$.

\section{Raising and lowering operators}

Recall that in the Lie algebra ${\rm so}(n,{\Bbb C})$ there
exist root elements
$E_{\alpha _1},E_{\alpha _2},\cdots ,E_{\alpha _k}$, corresponding to
simple roots, and root elements
$F_{\alpha _1},F_{\alpha _2},\cdots ,F_{\alpha _k}$, corresponding to
simple roots taken with sign minus. If $T'$ is a finite dimensional
irreducible representation of ${\rm so}(n,{\Bbb C})$ and $|{\bf m}\rangle$
is its weight vector, then
\[
T'(E_{\alpha _i}) |{\bf m}\rangle =\beta_{\bf m}
|{\bf m}+\alpha _i\rangle ,\ \ \ \ \
T'(F_{\alpha _i}) |{\bf m}\rangle =\gamma_{\bf m}
|{\bf m}-\alpha _i\rangle ,
\]
where $\beta_{\bf m}$ and $\gamma_{\bf m}$ are complex numbers.
In the algebra $U'_q({\rm so}_n)$ there exist no elements
similas to $E_{\alpha _j}$ and $F_{\alpha _j}$.
However, in finite dimensional representations of $U'_q({\rm so}_n)$
there exist operators having properties of the operators
$T'(E_{\alpha _i})$ and $T'(F_{\alpha _i})$. These operators
depend on a weight on which they act and are called {\it raising and
lowering operators} of the representation.
They are described as follows.

Let $T$ be a finite dimensional irreducible
representation of $U'_q({\rm so}_n)$ of the classical type and let
$|{\bf m}\rangle$ be its weight vector. If $n=2k$ we create the operators
$$
R_{\alpha _i}^{\bf m}=-T(I_{2i+2,2i-1})+q^{-(m_i+m_{i+1})/2}T(I_{2i+1,2i})
-{\rm i}q^{-m_i+1/2}T(I_{2i+2,2i}),\ \ \ \
$$    $$
\qquad\qquad -{\rm i}q^{-m_{i+1}-1/2}T(I_{2i+1,2i-1}) ,
\ \ \ \ i=1,2,\cdots ,k-1,   \eqno (14)
$$    $$
L_{\alpha _i}^{\bf m}=-T(I_{2i+2,2i-1})+q^{(m_i+m_{i+1})/2}T(I_{2i+1,2i})
+{\rm i}q^{m_i+1/2}T(I_{2i+2,2i})
$$   $$
\qquad\qquad
+{\rm i}q^{m_{i+1}-1/2}T(I_{2i+1,2i-1}),\ \ \ \ i=1,2,\cdots ,k-1,  \eqno
(15)
$$
and the operators
$$
R_{\alpha _k}^{\bf m}=T(I_{2k,2k-3})+q^{(-m_{k-1}+m_{k})/2}T(I_{2k-1,2k-2})
+{\rm i}q^{-m_{k-1}+1/2}T(I_{2k,2k-2})
$$    $$
\qquad\qquad
-{\rm i}q^{m_{k}-1/2}T(I_{2k-1,2k-3}) ,  \eqno (16)
$$     $$
L_{\alpha _k}^{\bf m}=-T(I_{2k,2k-3})+q^{(m_{k-1}-m_{k})/2}T(I_{2k-1,2k-2})
+{\rm i}q^{m_{k-1}+1/2}T(I_{2k,2k-2})
$$     $$
\qquad\qquad
+{\rm i}q^{-m_{k}-1/2}T(I_{2k-1,2k-3}) .  \eqno (17)
$$
If $n=2k+1$, then we create the operators (14), (15) and the operators
$$
R_{\alpha _k}^{\bf m}=T(I_{2k+1,2k-1})+{\rm i}q^{-m_k+1/2}T(I_{2k+1,2k}),
 \eqno (18)
$$    $$
L_{\alpha _k}^{\bf m}=T(I_{2k+1,2k-1})-{\rm i}q^{m_k+1/2}T(I_{2k+1,2k}).
 \eqno (19)
$$

If $T$ is a finite dimensional
representation of $U'_q({\rm so}_n)$ of the nonclassical type and
$|{\bf m}\rangle$ is its weight vector, then we create the operators
$$
R_{\alpha _i}^{\bf m}=-T(I_{2i+2,2i-1})+q^{-(m_i+m_{i+1})/2}T(I_{2i+1,2i})
-q^{-m_i+1/2}T(I_{2i+2,2i})
$$    $$
\qquad\qquad
-q^{-m_{i+1}-1/2}T(I_{2i+1,2i-1}),\ \ \ \ i=1,2,\cdots ,k-1,   \eqno (20)
$$    $$
L_{\alpha _i}^{\bf m}=-T(I_{2i+2,2i-1})+q^{(m_i+m_{i+1})/2}T(I_{2i+1,2i})
-q^{m_i+1/2}T(I_{2i+2,2i})
$$     $$
\qquad\qquad
-q^{m_{i+1}-1/2}T(I_{2i+1,2i-1}) , \ \ \ \ i=1,2,\cdots ,k-1 ,  \eqno (21)
$$
and the operators
$$
R_{\alpha _k}^{\bf m}=T(I_{2k,2k-3})+q^{(-m_{k-1}+m_{k})/2}T(I_{2k-1,2k-2})
+q^{-m_{k-1}+1/2}T(I_{2k,2k-2})
$$     $$
\qquad\qquad
+q^{m_{k}-1/2}T(I_{2k-1,2k-3}) ,   \eqno (22)
$$    $$
L_{\alpha _k}^{\bf m}=-T(I_{2k,2k-3})-q^{(m_{k-1}-m_{k})/2}T(I_{2k-1,2k-2})
+q^{m_{k-1}+1/2}T(I_{2k,2k-2})
$$     $$
\qquad\qquad
+q^{-m_{k}-1/2}T(I_{2k-1,2k-3})    \eqno (23)
$$
if $n=2k$. If $n=2k+1$, then we create the operators (20), (21) and
the operators
$$
R_{\alpha _k}^{\bf m}=T(I_{2k+1,2k-1})+q^{-m_k+1/2}T(I_{2k+1,2k}),
 \eqno (24)
$$    $$
L_{\alpha _k}^{\bf m}=T(I_{2k+1,2k-1})-q^{m_k+1/2}T(I_{2k+1,2k}).
 \eqno (25)
$$

The operators $R_{\alpha _k}^{\bf m}$ and $L_{\alpha _k}^{\bf m}$
correspond to the operators $T'(E_{\alpha _i})$ and $T'(F_{\alpha _i})$
of a representation $T'$ of the Lie algebra ${\rm so}(n,{\Bbb C})$,
respectively. We have
$$
R_{\alpha _k}^{\bf m}|{\bf m}\rangle =\beta _i
|{\bf m}+\alpha _i\rangle ,\ \ \ \ \
L_{\alpha _k}^{\bf m}|{\bf m}\rangle =\gamma _i
|{\bf m}-\alpha _i\rangle ,   \eqno (26)
$$
where $\alpha _i$ and $\gamma _i$ are complex numbers, which depend on
the representation of $U'_q({\rm so}_n)$.
Note that the relations (26) is not correct if we replace
the vector $|{\bf m}\rangle $ by some other weight vector
$|{\bf m}'\rangle$, since in such a case on the right hand side
we shall obtain, except of the vectors
$|{\bf m}'+\alpha _i\rangle  $ and $|{\bf m}'-\alpha _i\rangle$,
other weight vectors.

Formulas (14)--(17) and (20)--(23) for raising and lowering operators
follows from formulas of section 8 of the paper [8] if to restrict
the representation $T$ of $U'_q({\rm so}_n)$ to the subalgebras
$U'_q({\rm so}_4)$ generated by the elements $I_{2j,2j-1},
I_{2j+1,2j},I_{2j+2,2j+1}$, $j=1,2,\cdots ,k-1$.

Formulas (18), (19), (24) and (25) for raising and lowering operators
follows from formulas for raising and lowering operators for
irreducible representations of the algebra $U'_q({\rm so}_3)$
of the paper [8] if to restrict
the representation $T$ to the subalgebra
$U'_q({\rm so}_3)$ generated by the elements $I_{2k+1,2k}$ and
$I_{2k,2k-1}$.
\bigskip

{\bf Definition 3.} If $T$ is a finite dimensional irreducible
representation of the
algebra $U'_q({\rm so}_n)$, then a weight ${\bf m}$ of this representation
is called a {\it highest weight} if $R^{\bf m}_{\alpha _i}
|{\bf m}\rangle =0$, $i=1,2,\cdots ,k$. The corresponding vector
$|{\bf m}\rangle $ is called a {\it highest weight vector}.
\bigskip

Let us give a form of highest weights of irreducible representations of
the classical and of the nonclassical types. In order to determine
such a form we restrict the corresponding irreducible representations of
$U'_q({\rm so}_n)$ to the subalgebras $U'_q({\rm so}_4)$ and
$U'_q({\rm so}_3)$ and use the results of the papers [8] and [23].
As a result, we find that if a weight ${\bf m}\equiv (m_1,m_2,\cdots ,
m_k)$ of an irreducible representation $T$ of the classical type
is a highest weight, then the numbers $m_j$ are all integral or all
half-integral (but not integral) and satisfy the conditions
\[
m_1\ge m_2\ge \cdots \ge m_k\ \ \ \ \ {\rm if}\ \ \ \ \
n=2k+1
\]
and the conditions
\[
m_1\ge m_2\ge \cdots \ge m_{k-1}\ge | m_k|\ \ \ \ \ {\rm if}\ \ \ \ \
n=2k.
\]
The set of these highest weights coincides with the set of
highest weights of irreducible finite dimensional representations
of the Lie algebra ${\rm so}(n,{\Bbb C})$. These highest weights will be
calles {\it highest weights of the classical type}.

If a weight ${\bf m}\equiv (m_1,m_2,\cdots ,
m_k)$ of an irreducible representation $T$ of the nonclassical type
is a highest weight, then the numbers $m_j$ are all
half-integral (but not integral). In order to formulate the
classification theorem for representations of the nonclassical type we
shall need only highest weights ${\bf m}$ for which all $m_j$ are
positive. Such highest weights must satisfy
the conditions
\[
m_1\ge m_2\ge \cdots \ge m_k\ge 1/2.
\]
These highest weights will be called {\it highest weights of the
nonclassical type}.

It is well known that the root elements $E_{\alpha _i}$ and
$F_{\alpha _i}$ of the Lie algebra ${\rm so}(n,{\Bbb C})$
satisfy the relations
\[
[E_{\alpha _i}, F_{\alpha _i}]=2H_{\alpha _i},\ \ \ \ \
[E_{\alpha _i}, F_{\alpha _j}]=0,\ \ \ i\ne j.
\]
Instead of these relations for raising and lowering operators
of representations of the classical and nonclassical type
of the algebra $U'_q({\rm so}_{2k})$ we have the
relations
$$
(R_{\alpha _i}^{{\bf m}-\alpha _i}L_{\alpha _i}^{\bf m} -
L_{\alpha _i}^{{\bf m}+\alpha _i}R_{\alpha _i}^{\bf m} )
|{\bf m}\rangle =[2l]_q\{ (q-q^{-1})^2C_4
-(q^{2l}+q^{-2l})(q^2-q^{-2})\} |{\bf m}\rangle ,   \eqno (27)
$$   $$
(R_{\alpha _i}^{{\bf m}-\alpha _j}L_{\alpha _j}^{\bf m} -
L_{\alpha _j}^{{\bf m}+\alpha _i}R_{\alpha _i}^{\bf m} )
|{\bf m}\rangle =0,     \eqno (28)
$$
where $l=(m_i-m_{i+1})/2$ if $i\ne k$ and $l=(m_i+m_{i+1})/2$ if
$i=k$ and $C_4$ is the Casimir operator of the
subalgebra $U'_q({\rm so}_4)$ generated by the elements
$I_{2i,2i-1},I_{2i+1,2i},I_{2i+2,2i+1}$, which is given as
\[
C_4=q^{-1}I_{2i,2i-1}I_{2i+2,2i+1}-I_{2i+1,2i-1}I_{2i+2,2i}
+qI_{2i+1,2i}I_{2i+2,2i-1}.
\]

For the algebra $U'_q({\rm so}_{2k+1})$ we have the relations (27) for
$i\ne k$, (28) for $i\ne j$ and the relations
$$
(R_{\alpha _k}^{{\bf m}-\alpha _k}L_{\alpha _k}^{\bf m} -
L_{\alpha _k}^{{\bf m}+\alpha _k}R_{\alpha _k}^{\bf m} )
|{\bf m}\rangle =q[m_k]_q[m_k]_+(q-q^{-1})^2 |{\bf m}\rangle .  \eqno (29)
$$

\section{Classification theorems}

For finite dimensional irreducible representations of the classical type
the following theorem is true.
\bigskip

{\bf Theorem 2.}
(a) {\it Each irreducible finite dimensional representation of the
classical type
has a highest weight. A highest weight is unique (up to a constant).}

(b) {\it Irreducible finite dimensional representations with different
highest weights are not equivalent. Conversely, nonequivalent irreducible
finite dimensional representations of $U'_q({\rm so}_n)$ have
different highest weights.}
\bigskip

Existing of a highest weight is proved in the same way as in the
case of irreducible representations of the Lie algebra ${\rm so}(n, {\Bbb
C})$
by using Propositions 1 and 2. A proof of uniqueness of highest weight is
not simple. The relations (27)--(29) are used in this proof.

The assertion (b) is proved by using a proof of the similar assertion
for irreducible representations of the algebra $U'_q({\rm so}_4)$
from paper [8]. Namely, if $| {\bf m}\rangle$ is a highest weight
vector, then we act upon $| {\bf m}\rangle$ successively by the
corresponding operators $L^{{\bf m}'}_{\alpha _i}$, $i=1,2,\cdots ,k$.
Then, as in [8], we can find how the operators
$R^{{\bf m}'}_{\alpha _i}$ act upon weight vectors $| {\bf m}'\rangle$.
Therefore, by the method of the paper [8] we evaluate uniquely how the
operators $T(I_{2i+2,2i-1}),T(I_{2i+1,2i}),T(I_{2i+2,2i}), T(I_{2i+2,2i+1})$
act upon the corresponding weight vectors. Thus, a highest weight
determines uniquely (up to equivalence) the operators
$T(I_{j,j-1})$, $j=2,3,\cdots ,n$.
\medskip

Thus, in order to obtain a classification of irreducible finite
dimensional representations of the classical type of the algebra
$U'_q({\rm so}_n)$
we have to determine for which highest weights, described in
the previous section, there correspond such irreducible representations
with these highest weights.

It can be proved that the irreducible representation $T_{\bf m}$
of $U'_q({\rm so}_n)$ from the paper [5] are of the classical type
and has highest weight ${\bf m}$.
If we take all these irreducible representations $T_{\bf m}$,
then they give all highest weights ${\bf m}$, described in previous section
for irreducible representations of the classical type.
That is, for each highest weight ${\bf m}$ of the classical type
from the previous section there corresponds an irreducible representation
of $U'_q({\rm so}_n)$. Thus, we obtain the following classification
of irreducible representations of the classical type.
\bigskip

{\bf Theorem 3.} {\it
Irreducible finite dimensional representations of the classical type of the
algebra $U'_q({\rm so}_n)$ are in one-to-one correspondence with
highest weights of the classical type, described in the previous section.}
\bigskip

Thus,
irreducible finite dimensional representations of the classical type of the
algebra $U'_q({\rm so}_n)$ are in one-to-one correspondence with
irreducible finite dimensional representations of the Lie
algebra ${\rm so}_n$. The corresponding irreducible
representations of $U'_q({\rm so}_n)$ and of ${\rm so}_n$ act on
the same vector space. Moreover, when $q\to 1$, then operators
of an irreducible representation of $U'_q({\rm so}_n)$ tend to
the corresponding operators of the corresponding irreducible
representation of ${\rm so}_n$. This is a reason why the representations
of Theorem 3 are called representations of the classical type.

An analogue of Theorem 2 for irreducible representations of the nonclassical
type is formulated as follows.
\bigskip

{\bf Theorem 4.}
(a) {\it Each irreducible finite dimensional representation of the
nonclassical
type has a highest weight. A highest weight is unique (up to a constant).}

(b) {\it Irreducible finite dimensional representations with different
highest weights are not equivalent.}
\bigskip

This theorem is proved in the same way as Theorem 2.

In order to formulate the classification theorem for irreducible
representations of the nonclassical type we first formulatre the
following proposition.
\bigskip

{\bf Proposition 5.}
{\it If $T$ is an irreducible representation of the nonclassical type
and $G$ is the automorphism group of $U'_q({\rm so}_n)$  from
section 2, then the composition $T^{(g)}:=T\circ g$, $g\in G$, $g\ne e$,
is a representation of
the nonclassical type which is not equivalent to $T$.}
\bigskip

This proposition is proved by showing that spectrum of the operator
$T(I_{2i,2i-1})$ $(i=1,2,\cdots$, $k)$ coincides with the set
$[\frac 12]_+$, $[\frac 32]_+, \cdots ,[\frac s2]_+$ or with the set
$-[\frac 12]_+$, $-[\frac 32]_+, \cdots ,-[\frac s2]_+$, where
$s$ is some positive integer. In order to show this we use method of
mathematical induction. For $U'_q({\rm so}_4)$ this assertion is true
(see [8]). The induction is proved by using Wigner--Eckart theorem for
irreducible representations of the nonclassical type derived by
N. Iorgov (this theorem will be published).
\medskip

Thus, with every irreducible representation $T$ of the nonclassical type
we associate a set of irreducible representations $\{ T^{(g)}\ |\
g\in G\}$, consisting of $2^{n-1}$ pairwise nonequivalent
irreducible representations of the nonclassical type. In this set there
exists exactly one irreducible representation with highest weight
${\bf m}=(m_1, m_2,\cdots ,m_k)$ such that
$m_1\ge m_2\ge \cdots \ge m_k\ge \frac 12$.

For every highest weight of the nonclassical type ${\bf m}$ with
$m_1\ge m_2\ge \cdots \ge m_k\ge \frac 12$ there
exists constructed an irreducible representation of the nonclassical type
having ${\bf m}$ as its highest weight. Therefore, from above
reasoning we derive the following classification of irreducible
representations of the nonclassical type.
\bigskip

{\bf Theorem 5.} {\it
Irreducible representations of the nonclassical type of the algebra
$U'_q({\rm so}_n)$ are in one-to-one correspondence with pairs
$({\bf m},g)$, where ${\bf m}$ is a highest weight of the nonclassical
type with $m_1\ge m_2\ge \cdots \ge m_k\ge \frac 12$
and $g$ is an element of the automorphism group $G$.}
\bigskip

Note that irreducible representations of the nonclassical type have no
classical analogue. Namely, operators of representations of the
nonclassical type are singular at the point $q=1$.

\section{Irreducible representations of $U'_q({\rm so}_3)$}

This and the next sections are devoted to examples of the theory
described above. In this section we describe irreducible finite
dimensional representations of the algebra $U'_q({\rm so}_3$).

Irreducible finite dimensional representations of the classical type of
this algebra are given by nonnegative integral or half-integral number $l$.
The irreducible representation $T_l$, given by such a number $l$, acts
on $(2l+1)$-dimensional vector space ${\cal H}_l$ with a basis
$|l,m\rangle$, $m=-l,-l+1,\cdots ,l$. The operators $T_l(I_{21})$ and
$T_l(I_{32})$ are given by the formulas
$$
T_l(I_{21}) |l,m\rangle ={\rm i}[m]_q|l,m\rangle ,
$$     $$
T_l(I_{32}) |l,m\rangle ={\frac{1}{q^m+q^{-m}}}\left( [l-m]_q|l,m+1\rangle
-[l+m]_q|l,m-1\rangle \right) ,
$$
where $[a]_q$ denotes a $q$-number. Note that for these representations
we have
\[
{\rm Tr}\ T_l(I_{21})=0,\ \ \ \ \ \ {\rm Tr}\ T_l(I_{32})=0.
\]

Irreducible representations $T^{\epsilon _1, \epsilon _2}_n$ of the
nonclassical type are given by the numbers $\epsilon _i=\pm 1$
(they determine elements of the automorphism group $G$) and
by the integer $n=1,2,\cdots $. (According to section 6, these
representations are given by half-integral number $l$, but we replaced
$l$ by $n=l+1/2$.) The representation $T^{\epsilon _1, \epsilon _2}_n$
acts on $n$-dimensional vector space with the basis $|k\rangle$,
$k=1,2,\cdots ,n$. The operators $T^{\epsilon _1, \epsilon _2}_n(I_{21})$
and $T^{\epsilon _1, \epsilon _2}_n(I_{32})$ are given by the formulas
$$
T^{\epsilon _1, \epsilon _2}_n(I_{21}) |k\rangle =\epsilon _1
\frac{q^{k-1/2}+q^{-k+1/2}}{q-q^{-1}}|k\rangle ,
$$   $$
T^{\epsilon _1, \epsilon _2}_n(I_{32}) |1\rangle =
{\frac{1}{q^{1/2}-q^{-1/2}}}\left( \epsilon _2 [n]_q|1\rangle
+{\rm i}[n-1]_q |2\rangle \right) ,
$$    $$
T^{\epsilon _1, \epsilon _2}_n(I_{32}) |k\rangle =
\frac{1}{q^{k-1/2}-q^{-k+1/2}}
\left( {\rm i}[n-k]_q |k+1\rangle ) +
+{\rm i}[n+k-1]_q |k-1\rangle \right) ,
$$

These representations have the properties
\[
{\rm Tr}\ T^{\epsilon _1, \epsilon _2}_n(I_{21})\ne 0,\ \ \ \ \ \
{\rm Tr}\ T^{\epsilon _1, \epsilon _2}_n(I_{32})\ne 0.
\]
There exist 4 one-dimensional irreducible representations of the
nonclassical type. They are equivalent to
$T^{\epsilon _1, \epsilon _2}_1$, $\epsilon _i=\pm 1$.

Note that a proof of the fact that these representations of
$U'_q({\rm so}_3)$ exhaust all irreducible representations of this
algebra is given in [23].

\section{Irreducible representations of $U'_q({\rm so}_4)$}

Irreducible finite dimensional representations of the classical type
of the algebra
$U'_q({\rm so}_4)$ are given by two integral or two half-integral (but not
integral) numbers $r$ and $s$ such that $r\ge |s|$. These numbers
constitute the highest weight of the representation. We define the numbers
$j=(r+s)/2$ and $j'=(r-s)/2$ and denote the representation by $T_{jj'}$.
This representation acts on the vector space with the basis
\[
|k,l\rangle ,\ \ \ \ \ k=-j,-j+1,\cdots ,j,\ \ \  l=-j',-j'+1,\cdots ,j'.
\]
The operators $T_{jj'}(I_{i,i-1})$, $i=2,3,4$, act upon these vectors
by the formulas
$$
T_{jj'}(I_{21})|k,l\rangle ={\rm i}[k+l]_q|k,l\rangle , \ \ \ \
T_{jj'}(I_{43})|k,l\rangle ={\rm i}[k-l]_q|k,l\rangle ,
$$     $$
T_{jj'}(I_{32})|k,l\rangle =\frac{1}{(q^{k+l}+q^{-k-l})(q^{k-l}+q^{-k+l})}
\times
$$     $$
\qquad\qquad \times  \{ -(q^{j-l}+q^{-j+l}) [j'-l]_q |k,l+1\rangle +
(q^{j+l}+q^{-j-l}) [j'+l]_q |k,l-1\rangle +
$$    $$
\qquad\qquad
+(q^{j'-k}+q^{-j'+k}) [j-k]_q |k+1,l\rangle -
(q^{j'+k}+q^{-j'-k}) [j+k]_q |k-1,l\rangle \} .
$$

Irreducible finite dimensional representations of the nonclassical type
of the algebra
$U'_q({\rm so}_4)$ are given by two half-integral (but not
integral) numbers $r, s$ such that $r\ge s> 0$ and by the
numbers $\epsilon _1, \epsilon _2,\epsilon _3$, $\epsilon _i=\pm 1$,
which determine elements of the automorphism group $G$.
The numbers $r$ and $s$ constitute a highest weight of the
representation if $\epsilon _1=\epsilon _2=\epsilon _3=1$.
We define the numbers
$j=(r+s)/2$ and $j'=(r-s)/2$ and denote the corresponding representations by
$T^{\epsilon _1\epsilon _2,\epsilon _3}_{jj'}$.

If $(r,s)$ runs over all highest weights of the nonclassical type
with $r\ge s>0$, then $j$ and $j'$ run over the values
\[
\textstyle
j=0,1,2,\cdots ,\ \ j'={\frac 12},{\frac 32}, {\frac 53},\cdots
\ \ \ \ \ {\rm or}\ \ \ \ \
j={\frac 12},{\frac 32}, {\frac 53},\cdots ,\ \ j'=0,1,2,\cdots .
\]
The representation
$T^{\varepsilon _1,\varepsilon _2,\varepsilon _3}_{jj'}$
acts on the vector space ${\cal H}$ with the basis
\[
\textstyle
|k,l\rangle ,\ \ \ \ k=j,j-1,\cdots ,{\frac 12},\ \ \
l=j',j'-1,\cdots ,-j',
\]
if $j'$ is integral and with the basis
\[
\textstyle
|k,l\rangle ,\ \ \ \ k=j,j-1,\cdots ,-j,\ \ \
l=j',j'-1,\cdots ,{\frac 12},
\]
if $j$ is integral. The representations are given by the formulas
$$
T^{\varepsilon _1,\varepsilon _2,\varepsilon _3}_{jj'}(I_{21})
|k,l\rangle = \varepsilon _1[k+l]_+|k,l\rangle , \ \ \ \
T^{\varepsilon _1,\varepsilon _2,\varepsilon _3}_{jj'}(I_{43})
|k,l\rangle = \varepsilon _2 [k-l]_+ |k,l\rangle ,
$$    $$
T^{\varepsilon _1,\varepsilon _2,\varepsilon _3}_{jj'} (I_{32})
|k,l\rangle  =\frac{1}{[k+l]_q[k-l]_q(q-q^{-1})} \{
-{\rm i}[j'-l]_q[j-l]_q|k,l+1\rangle +
$$   $$
\qquad
+ {\rm i}[j'+l]_q[j+l]_q |k,l-1\rangle -{\rm i} [j'-k]_q[j-k]_q
|k+1,l\rangle + {\rm i} [j'+k]_q[j+k]_q |k-1,l\rangle \},
$$
where $k\ne {\frac 12}$ if $j$ is half-integral and
$l\ne {\frac 12}$ if $j'$ is half-integral, and by
$$
T^{\varepsilon _1,\varepsilon _2,\varepsilon _3}_{jj'}(I_{32})
| {\textstyle \frac 12} ,l\rangle  =\frac{1}
{ [l+\frac 12 ] [l-\frac 12 ] (q-q^{-1})}
\{ -{\rm i} [j-l]_q[j'-l]_q |{\textstyle \frac 12} ,l+1\rangle  +
$$    $$
\qquad
\textstyle
+{\rm i} [j+l]_q [j'+l]_q | \frac 12 ,l-1\rangle  -
{\rm i} [j'-\frac 12 ]_q[j-\frac 12 ]_q |\frac 32 ,l\rangle  +
{\rm i} [j'+\frac 12 ]_q [j+\frac 12 ]_q
\varepsilon _3 (-1)^l |\frac 12 ,-l\rangle \}
$$
if $j$ is half-integral and by
$$
T^{\varepsilon _1,\varepsilon _2,\varepsilon _3}_{jj'}(I_{32})
\textstyle
|k, \frac 12 \rangle  =\frac{1}
{ [k+\frac 12 ]_q [k-\frac 12 ]_q (q-q^{-1})}
\{ -{\rm i} [j-\frac 12 ]_q[j'-\frac 12 ]_q
|k, {\textstyle \frac 32} \rangle  +
{\rm i} [j{+}\frac 12 ]_q [j'{+}\frac 12 ]_q \times
$$    $$
\qquad
\textstyle
\times
\varepsilon _3 (-1)^k|-k, \frac 12 \rangle  -
{\rm i} [j'{-}k]_q[j{-}k]_q |k+1,\frac 12 \rangle  +
{\rm i} [j'{+}k]_q [j{+}k]+q |k-1,\frac 12 \rangle \}
$$
if $j'$ is half-integral.

Note that a proof of the fact that these representations of
$U'_q({\rm so}_4)$ exhaust all irreducible representations of this
algebra is given in [8].

\subsection*{Acknowledgement}

The research described in this paper was made possible in part by
Award No. UP1-2115 of the U.S. Civilian Research and Development
Foundation for Independent States of the Former Soviet Union (CRDF).


\begin{thebibliography}{99}

\bibitem{1}Jimbo M., A $q$-analogue of $U(g)$ and the Yang--Baxter
equation, {\em Lett. Math. Phys.}, 1985, V. 10, 63--69.
\bibitem{2}Drinfeld V. G., Hopf algebras and the quantum Yang--Baxter
equation, {\em Sov. Math. Dokl.}, 1985, V. 32, 354--258.
\bibitem{3}Jantzen J. C., {\it Lectures on Quantum Groups},
Aner. Math. Soc., Providence, RI, 1996.
\bibitem{4}Klimyk A. and Schm\"{u}dgen K., {\it Quantum Groups and Their
Representations}, Berlin, Springer, 1997.
\bibitem{5}Gavrilik A. M. and Klimyk A. U., $q$-Deformed orthogonal
and pseudo-orthogonal algebras and their representations,
{\em Lett. Math. Phys.}, 1991, V. 21, 215--220.
\bibitem{6}Fairlie D. B., Quantum deformations of $SU(2)$,
{\em J. Phys. A: Math. Gen.}, 1990, V. 23, L183--L187.
\bibitem{7}Odesskii A., An analogue of the Sklyanin algebra,
{\em Funct. Anal. Appl.}, 1986, V. 20, 152--154.
\bibitem{8}Havl\'{i}\v{c}ek M., Klimyk A. U. and Po\v{s}ta S.,
Representations of the $q$-deformed algebra $U'_q({\rm so}_4)$,
{\em J. Math. Phys.}, 2001, V. 42, N. 9.
\bibitem{9}Noumi M., Macdonald's symmetric polynomials as zonal
spherical functions on quantum homogeneous spaces, {\em Adv. Math.},
1996, V. 123, 16--77.
\bibitem{10}Iorgov N. Z.,
On tensor products of representations of the non-standard
$q$-deformed algebra $U'_q({\rm so}_n)$, {\em J. Phys. A: Math. Gen.},
2001, V. 34, 3095--3108.
\bibitem{11}Gavrilik A. M. and Iorgov N. Z., $q$-Deformed
algebras $U_q({\rm so}_n)$ and their representations,
{\em Methods of Funct. Anal. Topology}, 1997, V. 3, N. 4, 51--63.
\bibitem{12}Iorgov N. Z. and Klimyk A. U., Nonclassical type
representations of the $q$-deformed algebra $U'_q({\rm so}_n)$,
{\em Czech. J. Phys.}, 2000, V. 50, 85--90.
\bibitem{13}Havl\'{i}\v{c}ek M., Klimyk A. U. and Po\v{s}ta S.,
Representations of the cyclically symmetric $q$-deformed
algebra ${\rm so}_q(3)$, {\em J. Math. Phys.}, 1999, V. 40, 2135--2161.
\bibitem{14}Samoilenko Yu. S. and Turowska L.,
Semilinear relations and $*$-representations of deformations of
$SO(3)$, in {\it Quantum Groups and Quantum Spaces}, Banach Center
Publications, Vol. 40, Warsaw, 1997, p. 21--43.
\bibitem{15}Klimyk, A. U., Nonstandard $q$-deformation of the
universal enveloping algebra $U({\rm so}_n)$, in "{\it Quantum
Theory and Symmetry}", Editors H.-D. Doebner {\it et al}, Singapore,
World Scientific, 2000, 459--463.
\bibitem{16}Gavrilik A. M. and Iorgov N. Z.,
Representations of the nonstandard algebras $U_q({\rm so}(n))$ and
$U_q({\rm so}(n,1))$ in Gel'fand--Tsetlin basis, {\em Ukr. J. Phys.},
1998, V. 43, N. 7, 791--797.
\bibitem{17}Iorgov N. Z. and Klimyk A. U., The nonstandard
deformation $U'_q({\rm so}_n)$ for $q$ a root of unity,
{\em Methods of Funct. Anal. Topology}, 1999, V. 3, N. 6, 56--73.
\bibitem{18}Chekhov L. O.
and Fock V. V., Observables in 3D gravity and geodesic
algebras, {\em Czech. J. Phys.}, 2000, V. 50, N. 11, 1201--1208.
\bibitem{19}Nelson J. and Regge T., 2+1 quantum gravity,
{\em Phys. Lett.} B, 1991, V. 272, 213--216.
\bibitem{20}Nelson J. and Regge T., 2+1 gravity for genus $>1$,
{\em Commun. Math. Phys.}, 1991, V. 141, 211--223.
\bibitem{21}Bullock D. and Przytycky J. H., Multiplicative
structure of Kauffman bracket skein module quantization,
{\bf e}-print: math.QA/9902117.
\bibitem{22}Gavrilik A. M., The use of quantum algebras in quantum gravity,
{\em Proc. of Inst. of Mathem. of NAS of Ukraine}, 2000, V. 30, N. 2,
304--309.
\bibitem{23}Havl\'{i}\v{c}ek M. and Po\v{s}ta S., On the classification
of irreducible finite dimensional representations of $U'_q({\rm so}_3)$
algebra, {\em J. Math. Phys.}, 2001, V. 42, 472--491.




\end{thebibliography}
\end{document}